\documentclass[pdflatex,sn-mathphys-num]{sn-jnl}% Math and Physical Sciences Author Year Reference Style\usepackage{graphicx} % Required for inserting images

\usepackage{graphicx}%
\usepackage{multirow}%
\usepackage{amsmath,amssymb,amsfonts}%
\usepackage{amsthm}%
\usepackage{xcolor}%
\usepackage{algorithm}%
\usepackage{algpseudocodex}%
\usepackage{listings}%
\usepackage{hyperref}

\newcommand{\twovec}[2]{\begin{bmatrix}#1 \\ #2 \end{bmatrix}}

\newcommand{\infnorm}[1]{\left\| #1 \right\|_{\infty}}

\theoremstyle{thmstyleone}%
\newtheorem{theorem}{Theorem}%  meant for continuous numbers
\theoremstyle{thmstyleone}%

\theoremstyle{thmstylethree}%

\begin{document}

\title{Stable evaluation of derivatives for barycentric and continued fraction representations of rational functions}
\date{January 2026}

\author*[1]{\fnm{Tobin A.} \sur{Driscoll}}\email{driscoll@udel.edu}

\author[1]{\fnm{Yuxing} \sur{Zhou}}\email{zhouyx@udel.edu}
\equalcont{These authors contributed equally to this work.}

\affil[1]{\orgdiv{Department of Mathematical Sciences}, \orgname{University of Delaware}, \orgaddress{ \city{Newark}, \postcode{19716}, \state{DE}, \country{USA}}}

\abstract{Effective algorithms for approximation by rational functions exist for both barycentric and Thiele continued fraction (TCF) representations. We present the first numerically stable methods for derivative evaluation in the barycentric representation, including an $O(n)$ algorithm for all derivatives. We also extend an earlier $O(n)$ algorithm for evaluation of the TCF first derivative to higher orders. Numerical experiments confirm the robustness and efficiency of the proposed methods.}

\keywords{approximation, barycentric representation, continued fractions, numerical stability}

\maketitle

\section{Introduction}
\label{sec:intro}

Interest in computational approximation of functions by rational functions has grown since the advent of the AAA algorithm~\cite{NakatsukasaAAAAlgorithm2018}, which has found applications in model reduction, constructive approximation, PDEs, and more, as surveyed by Nakatsukasa and Trefethen~\cite{NakatsukasaApplicationsAAA2026}. The AAA algorithm computes an SVD to find weights for the barycentric interpolation representation of a rational function, thereby minimizing a linearized approximation error. This approach is embedded within a greedy iteration for selecting new interpolation nodes. More recently, a similar greedy approach to node selection using Thiele continued fraction (TCF) interpolation was suggested and demonstrated by Salazar Celis~\cite{SalazarCelisNumericalContinued2024} and shown to be 3--10 times faster than AAA to construct by Driscoll and Zhou~\cite{DriscollGreedyThiele2025pre}. Both AAA and greedy TCF produce interpolants on $n$ nodes that can be evaluated stably in $O(n)$ elementary operations.  

Derivatives of rational approximants are useful in many applications. Related work on derivatives of linear barycentric rational interpolants includes transformed Chebyshev rational interpolation~\cite{BaltenspergerExponentialConvergence1999} and convergence-rate results for derivatives of Floater--Hormann interpolants~\cite{BerrutConvergenceRates2011}. For barycentric rational interpolants, Schneider and Werner~\cite{SchneiderNewAspects1986} gave $O(n)$ formulas for the derivatives---one for a generic point, and a different one at nodes. However, the generic formula is numerically unstable near nodes~\cite{NakatsukasaApplicationsAAA2026}. For TCF interpolants, Driscoll and Zhou~\cite{DriscollGreedyThiele2025pre} gave an $O(n)$ iterative formula for the first derivative based on forward-mode automatic differentiation without investigation of its numerical stability.

In this work, we present a double-summation formula for the first derivative of the barycentric representation that removes the dominant cancellation error but requires $O(n^2)$ operations. We then derive a modified $O(n)$ algorithm for all derivatives that avoids the dominant cancellation. For TCF interpolants, we show that the $O(n)$ iterative formula given by Driscoll and Zhou~\cite{DriscollGreedyThiele2025pre} is stable in practice and straightforward to generalize to higher derivatives. 

The rest of this paper is organized as follows. In \autoref{sec:bary}, we review the barycentric representation and the original formula of Schneider and Werner for its derivative. We point out the source of subtractive cancellation and then present methods that avoid it. In \autoref{sec:cfrac}, we review the Thiele continued fraction representation and the iteration of Driscoll and Zhou to compute its first derivative, which we then generalize to higher-order derivatives. Numerical experiments are reported in \autoref{sec:experiments}, and conclusions are drawn in \autoref{sec:conclusions}.

\section{Barycentric representation}
\label{sec:bary}

The barycentric representation of a rational function with $n$ distinct nodes $z_1,z_2,\ldots,z_{n}$ is
\begin{equation}
r(z) = \frac{\displaystyle \sum_{k=1}^{n} \frac{w_k f_k}{z-z_k}}{\displaystyle  \sum_{k=1}^{n} \frac{w_k}{z-z_k}} = \frac{N(z)}{D(z)},
\label{eq:bary}
\end{equation}
where $f_k = f(z_k)$ are interpolated function values at the nodes and $w_k$ are the barycentric weights, which can be used to improve approximation properties.

\begin{figure}
    \centering
    \includegraphics[width=0.9\textwidth]{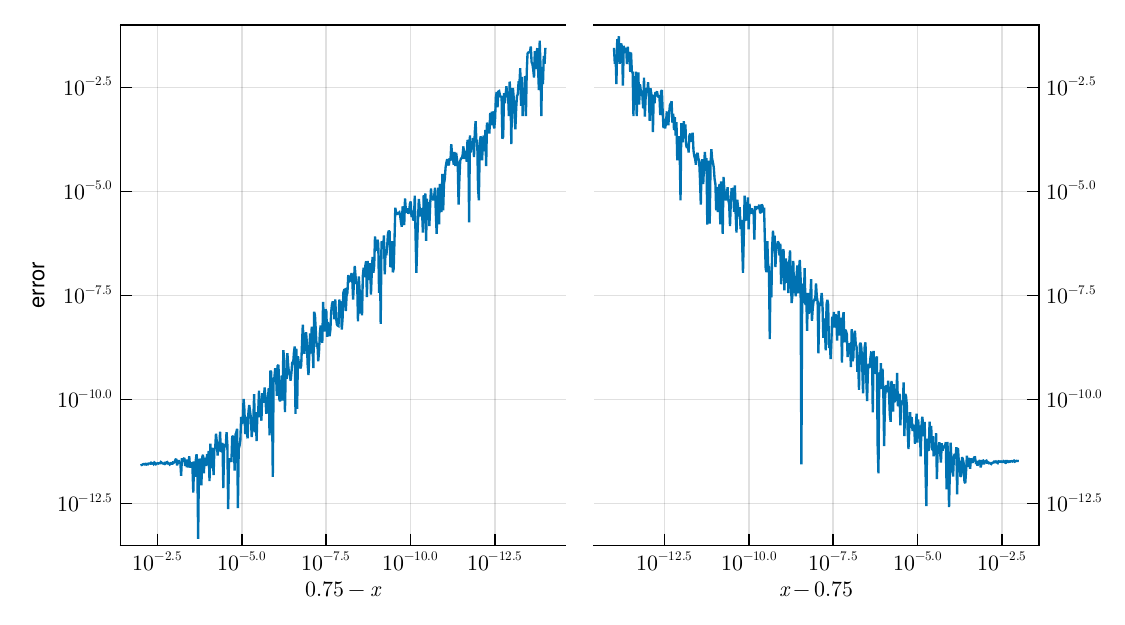}\\
    \includegraphics[width=0.9\textwidth]{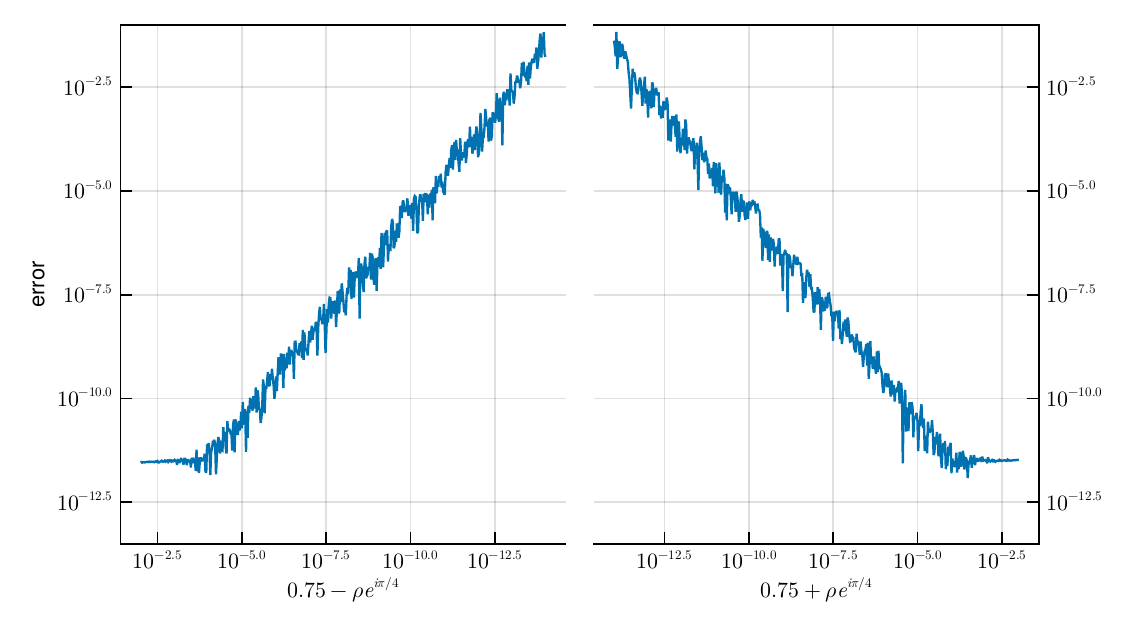}%
    \caption{Absolute error near a node using~\eqref{eq:sw1} to evaluate the first derivative of an AAA rational approximation $r(x)$ to $f(x) = e^x$ on $[-1,1]$. The approximation to $f$ is of type $(5,5)$ and accurate to within $5\times 10^{-13}$. Top: Errors along the real axis near the evaluation point $x=0.75$. Bottom: Errors along the line passing through $x=0.75$ at an angle of $\pi/4$ with the real axis.}
    \label{fig:swderiv_error}
\end{figure}

Implicit differentiation of $r(z)D(z) = N(z)$ leads to 
\begin{equation}
    r'(z) = \frac{N'(z)- r(z)D'(z)}{D(z)}.
    \label{eq:sw1}
\end{equation}
This formula was suggested by Schneider and Werner~\cite{SchneiderNewAspects1986} for computing the first derivative of a barycentric representation at points other than nodes. However, it is numerically unstable when $z$ is close to one of the nodes, say $z_j$. Define $\epsilon=z-z_j$ and
\begin{equation}
    N_j(z) = \sum_{\substack{k=1 \\ k\neq j}}^n \frac{w_k f_k}{z-z_k}, \quad 
    D_j(z) = \sum_{\substack{k=1 \\ k\neq j}}^n \frac{w_k}{z-z_k}.
    \label{eq:N1D1}
\end{equation}
We expand
\begin{align}
    r(z) &= \left[ \frac{w_jf_j}{\epsilon} + N_j(z) \right] \left[ \frac{w_j}{\epsilon} + D_j(z) \right]^{-1}\\ 
    & = \left[ f_j + \frac{\epsilon}{w_j} N_j(z) \right] \left[1 +  \frac{\epsilon}{w_j} D_j(z) \right]^{-1}\\ 
    &= f_j + \frac{N_j - f_jD_j}{w_j} \, \epsilon + O(\epsilon^2).
\end{align}
Thus,
\begin{equation}
   N'(z)- r(z)D'(z) = \left[ \frac{-w_j f_j}{\epsilon^2} + N_j'(z)\right] -  
   \Bigl[ f_j + O(\epsilon) \Bigr] \left[ \frac{-w_j}{\epsilon^2} + D_j'(z) \right],
\end{equation}
which requires an unstable subtractive cancellation of the dominant terms. \autoref{fig:swderiv_error} illustrates the result of this subtractive cancellation when the rational function $r$ is an AAA approximation to $f(x) = e^x$ on $[-1,1]$ of type $(5,5)$, accurate to within $5\times 10^{-13}$. The error in evaluating $r'(x)$ using~\eqref{eq:sw1} grows as $x$ approaches the node at $0.75$ along two different lines in the complex plane, reaching more than $0.01$ at a distance of $10^{-14}$ from the node.\footnote{Errors along a line parallel to the imaginary axis (not shown in the figure) have a similar maximum envelope, but they oscillate rapidly due to an effect we believe is related to the conjugate symmetry of the situation.}

One way to avoid the cancellation instability is to multiply \eqref{eq:sw1} through by $D$ and rewrite the numerator as the double sum
\begin{equation}
    N'(z)D(z) - N(z)D'(z) = \sum_{j=1}^n \sum_{k=1}^n \frac{w_j w_k (f_k - f_j)}{(z-z_j)^2 (z-z_k)}.
    \label{eq:double-sum}
\end{equation}
Numerical experiments confirm that this formula resolves the instability, but since it has $O(n^2)$ complexity, we next suggest an alternative that is $O(n)$.

Consider the representation $r=\tilde{N}/\tilde{D}$, with
\begin{equation}
    \tilde{N}(z) = \epsilon N(z), \quad \tilde{D}(z) = \epsilon D(z),
\end{equation}
still using $\epsilon=z-z_j$. Then 
\begin{equation}
    r'(z) = \frac{\tilde{N}'(z)- r(z)\tilde{D}'(z)}{\tilde{D}(z)}.
    \label{eq:sw1mod}
\end{equation}
Expanding this numerator leads to
\begin{equation}
    \tilde{N}'(z)- r(z)\tilde{D}'(z) = \left[N_j(z) +  \epsilon N_j'(z)  \right] - \bigl[ f_j + O(\epsilon) \bigr] \bigl[ D_j(z) +  \epsilon D_j'(z)  \bigr],
\end{equation}
which has no intrinsic subtractive cancellation. Hence, we can use~\eqref{eq:sw1mod} stably instead of~\eqref{eq:sw1} at all points that are close to $z_j$. Note that as $z\to z_j$, we obtain
\begin{equation}
    r'(z_j) = \frac{N_j(z_j) - f_j D_j(z_j)}{w_j} = \frac{1}{w_j} \sum_{k=1, k\neq j}^n \frac{w_k(f_k - f_j)}{z_j - z_k},
\end{equation}
which is the Schneider--Werner alternate formula for the derivative at a node.

For higher-order derivatives, we can apply Leibniz's product rule formula to $r \tilde{D} = \tilde{N}$ and solve to obtain
\begin{equation}
    r^{(m)}(z) = \frac{1}{\tilde{D}(z)} \left[\tilde{N}^{(m)}(z) - \sum_{k=1}^{m} \binom{m}{k} r^{(m-k)}(z) \tilde{D}^{(k)}(z)\right], \quad m \ge 1,
    \label{eq:rderiv}
\end{equation}
where
\begin{equation}
    \tilde{N}^{(m)}(z) = \epsilon N^{(m)}(z) + m N^{(m-1)}(z), \quad \tilde{D}^{(m)}(z) = \epsilon D^{(m)}(z) + m D^{(m-1)}(z).
    \label{eq:tildederiv}
\end{equation}
In general, when evaluating the derivative at a point $z$, we can find the index $j$ of the node closest to $z$ in $O(n)$ time and use the modified formula~\eqref{eq:sw1mod} with $\epsilon = z - z_j$. An outline for evaluation of barycentric derivatives is given in \autoref{alg:baryderiv}.

\begin{algorithm}
    \caption{Stable evaluation of derivatives of barycentric rational functions}
    \label{alg:baryderiv}
\begin{algorithmic}[1]
    \Require Evaluation point $\zeta$, max derivative order $\mu$, and $n$-vectors of nodes $z$, function values $f$, and weights $w$ 
    \State $\delta_k \gets \zeta - z_k, \quad k=1,2,\ldots,n$
    \State $j \gets \arg \displaystyle \min_k |\delta_k|$
    \State $N_j^{(m)} \gets 0,\; D_j^{(m)} \gets 0, \quad m=0,1,\ldots,\mu$
    \For{$k=1,\ldots,n$, except $j$}
    \Comment{Calculate $N_j$ and $D_j$ from~\eqref{eq:N1D1}, and their derivatives }
        \State{$d \gets w_k / \delta_k$}
        \State{Add $f_k d$ to $N_j^{(0)}$}
        \State{Add $d$ to $D_j^{(0)}$}
        \For{$m=1,\ldots,\mu$}
            \State{$d \gets -m d / \delta_k$}
            \State{Add $f_k d$ to $N_j^{(m)}$}
            \State{Add $d$ to $D_j^{(m)}$}
        \EndFor
    \EndFor
    \State{$\tilde{N}^{(0)} \gets w_j f_j + \delta_j N_j^{(0)}$}
    \State{$\tilde{D}^{(0)} \gets w_j + \delta_j D_j^{(0)}$}
    \State{$r^{(0)} \gets \tilde{N}^{(0)} / \tilde{D}^{(0)}$}
    \For{$m=1,\ldots,\mu$}
    \Comment{Apply~\eqref{eq:rderiv} and~\eqref{eq:tildederiv}}
        \State{$\tilde{N}^{(m)} \gets m N_j^{(m-1)} + \delta_j N_j^{(m)}$} 
        \State{$\tilde{D}^{(m)} \gets m D_j^{(m-1)} + \delta_j D_j^{(m)}$}
        \State{$r^{(m)} \gets \left[ \tilde{N}^{(m)} - \sum_{j=1}^{m} \binom{m}{j} r^{(m-j)} \tilde{D}^{(j)} \right] / \tilde{D}^{(0)}$}
    \EndFor
    \State \Return{$r^{(0)}, r^{(1)}, \ldots, r^{(\mu)}$}
\end{algorithmic}
\end{algorithm}

\section{Continued fraction representation}
\label{sec:cfrac}

The Thiele continued fraction (TCF) representation of a rational function with $n$ distinct nodes $z_1,z_2,\ldots,z_{n}$ is
\begin{equation}
    r(z) = w_1 + \frac{z-z_1}{w_2 + \cfrac{z-z_2}{w_3 + \dots \cfrac{z-z_{n-1}}{w_n}}}.
    \label{eq:tcf}
\end{equation}
The TCF coefficients $w_k$ (which are unrelated to the weights in the barycentric formula) can be computed iteratively for $k=1,\ldots,n$ through
\begin{equation}
    \begin{aligned}
        t_1 &= f_k, \\
        t_{i+1} &= \frac{z_k - z_i}{-w_i + t_{i}}, \quad i=1,2,\ldots,k-1, \\
        w_k &= t_k,
    \end{aligned}
\end{equation}

Jones and Thron~\cite{jonesNumericalStability1974} recommended the tail-ordered evaluation of $r(z)$. In terms of the standard generic continued fraction expression
\begin{equation}
    r(z) = b_0 + \frac{a_1}{b_1 + \cfrac{a_2}{b_2 + \dots \cfrac{a_{n-2}}{b_{n-1} + \cfrac{a_{n-1}}{b_n}}}},
    \label{eq:cfgeneric}
\end{equation}
the evaluation procedure is 
\begin{equation}
    \begin{aligned}
        u_{n+1} &= 0, \\
        u_k &= \frac{a_k}{b_k + u_{k+1}}, \quad k=n,n-1,\ldots,1, \\
        r(z) &= u_0 = b_0 + u_1.
    \end{aligned}
    \label{eq:tcf-classic}
\end{equation}
Driscoll and Zhou, however, found~\cite{DriscollGreedyThiele2025pre} that an alternative iteration was faster in floating-point arithmetic, due to requiring just a single division:
\begin{equation}
    \begin{aligned}
        \twovec{p_{n+1}}{q_{n+1}} &= \twovec{1}{0}, \\
        \twovec{p_k}{q_k} &= \twovec{b_k p_{k+1} + q_{k+1}}{a_k p_{k+1}}, \quad k=n, n-1,\ldots,1,\\
        \twovec{p_0}{q_0} &= \twovec{b_0 p_{1} + q_{1}}{p_{1}}, \\
        r(z) &= \frac{p_0}{q_0}.
    \end{aligned}
    \label{eq:tcf-onediv}
\end{equation}
We now show that these alternatives have the same numerical stability.
\begin{theorem}
Using the standard model of floating-point arithmetic and stability~\cite{highamAccuracyStability2002}, the iterations~\eqref{eq:tcf-classic} and~\eqref{eq:tcf-onediv} have the same stability.
\end{theorem}
\begin{proof}
    First, we establish that $u_k = q_k/p_k$ in exact arithmetic for $k=1,\dots,n$. This is trivial for $k=n+1$. For $k=n,n-1,\ldots,1$, we have by induction that
    \begin{equation}
        \frac{q_k}{p_k} = \frac{a_k p_{k+1}}{b_k p_{k+1} + q_{k+1}} = \frac{a_k p_{k+1}}{b_k p_{k+1} + u_{k+1} p_{k+1}} = u_k.
        \label{eq:ratio}
    \end{equation}

    Multiplications and divisions in floating-point have condition number 1 and do not affect stability. By~\eqref{eq:ratio}, the addition step $b_k p_{k+1} + q_{k+1}$ in~\eqref{eq:tcf-onediv} is a multiple $p_{k+1}$ of the addition $b_k + u_{k+1}$ in~\eqref{eq:tcf-classic}, so that step has the same relative condition number in both algorithms.
\end{proof}

Both~\eqref{eq:tcf-classic} and~\eqref{eq:tcf-onediv} are easily differentiated using forward-mode automatic differentiation. Driscoll and Zhou~\cite{DriscollGreedyThiele2025pre} gave the algorithm based on~\eqref{eq:tcf-onediv} for the first derivative, which we extend in \autoref{alg:tcfderiv} to higher derivatives. Defining $b_k=w_{k+1}$ and $a_k = z - z_k$, and stopping one term earlier in~\eqref{eq:cfgeneric}, we note that the $m$th derivatives of $p_k$ and $q_k$ satisfy the recurrences
\begin{equation}
\begin{split}
    p_k^{(m)} &= w_{k+1} p_{k+1}^{(m)} + q_{k+1}^{(m)}, \\
    q_k^{(m)} &= (z - z_k) p_{k+1}^{(m)} + m p_{k+1}^{(m-1)},
    \label{eq:tcf-deriv}
\end{split}
\end{equation}
initialized by $p_{n}^{(m)} = q_{n}^{(m)} = 0$ for $m >0$, and $p_{n}^{(0)} = 1$, $q_{n}^{(0)} = 0$. For the final iteration with $k=0$, the first equation is the same, and $q_0^{(m)} = p_{1}^{(m)}$. The derivatives of $r$ then follow as in~\eqref{eq:rderiv}, with $p$ and $q$ replacing $\tilde{N}$ and $\tilde{D}$.

In practice, we observe that it is sometimes important to guard against numerical underflow when applying~\eqref{eq:tcf-deriv}, by rescaling $p_k^{(m)}$ and $q_k^{(m)}$ by a common factor if they are both very small. This does not affect the final evaluation of $r^{(m)}$, since the scaling factors cancel in the ratio.

\begin{algorithm}
\caption{Evaluation of a Thiele continued fraction and its derivatives.}
\label{alg:onediv}
\begin{algorithmic}
    \Require Evaluation point $\zeta$, max derivative order $\mu$, and $n$-vectors of nodes $z$ and weights $w$ 
    \State{$p_n^{(0)}\gets 1$, $q_n^{(0)}\gets 0$}
    \State{$p_n^{(m)}\gets 0$, $q_n^{(m)}\gets 0$ for $m=1,\ldots,\mu$}
    \For{$k=n-1,\ldots,1$}
        \State{$p_k^{(0)} \gets w_{k+1} p_{k+1}^{(0)} + q_{k+1}^{(0)}$}
        \State{$q_k^{(0)} \gets (\zeta - z_k) p_{k+1}^{(0)}$}
        \For{$m=1,\ldots,\mu$}
            \State{$p_k^{(m)} \gets w_{k+1} p_{k+1}^{(m)} + q_{k+1}^{(m)}$}
            \State{$q_k^{(m)} \gets (\zeta - z_k) p_{k+1}^{(m)} + m p_{k+1}^{(m-1)}$}
        \EndFor
    \EndFor
    \State{$p_0^{(0)} \gets w_{1} p_{1}^{(0)} + q_{1}^{(0)}$, $q_0^{(0)} \gets p_{1}^{(0)}$}
    \State{$r^{(0)} \gets p_0^{(0)} / q_0^{(0)}$}
    \For{$m=1,\ldots,\mu$}
        \State{$p_0^{(m)} \gets w_{1} p_{1}^{(m)} + q_{1}^{(m)}$}
        \State{$q_0^{(m)} \gets p_{1}^{(m)}$}
        \State{$r^{(m)} \gets \left[ p_0^{(m)} - \sum_{j=1}^{m} \binom{m}{j} r^{(m-j)} q_0^{(j)} \right] / q_0^{(0)}$}
    \EndFor
    \State \Return{$r^{(0)}, r^{(1)}, \ldots, r^{(\mu)}$}
\end{algorithmic}
\label{alg:tcfderiv}
\end{algorithm}

\section{Numerical experiments}
\label{sec:experiments}

We have implemented \autoref{alg:baryderiv} and \autoref{alg:tcfderiv} in version 0.3.8 of the RationalFunctionApproximation.jl package~\cite{driscollRationalFunctionApproximationjlRational2023} for the Julia programming language. Here, we report the results of tests of these implementations on the functions listed in \autoref{tab:testfun-interval} over the interval $[-1,1]$ and the functions listed in \autoref{tab:testfun-circle} on the unit circle. Each function $f$ is scaled so the derivative being measured has max-norm 1 over the domain of interest. 

\begin{table}
    \centering
    \begin{tabular}{ll}
        Function & Characteristics\\
        \hline
        $f_E(x) = c\cdot \exp(\sin(x))$ & Entire function \\
        $f_C(x) = c\cdot\cos(20x)$ & Highly oscillatory \\
        $f_T(x) = c\cdot\tanh(x/\epsilon)$ & Meromorphic with nearest poles at $\pm i \pi \epsilon/2$ \\
        $f_L(x) = c\cdot\log(1 + \epsilon - x)$ & Branch point at $1 + \epsilon$ \\
        $f_A(x) = c\cdot\arctan(x/\epsilon)$ & Branch points at $\pm i \epsilon$ 
    \end{tabular}
    \caption{Test functions used in numerical experiments on the interval $[-1,1]$. The scaling constant $c$ is always chosen so that $\infnorm{f^{(m)}} = 1$ over the domain, where $m$ is the derivative order being tested.}
    \label{tab:testfun-interval}
\end{table}

\begin{table}
    \centering
    \begin{tabular}{ll}
        Function & Characteristics\\
        \hline
        $f_E(z) = c\cdot \exp(\sin(z))$ & Entire function \\
        $f_O(z) = c\cdot\cos(z^{10})$ & Highly oscillatory \\
        $f_L(z) = c\cdot\log(1 + \epsilon - z)$ & Branch point at $1 + \epsilon$ \\
        $f_S(z) = c\cdot\sqrt{1 - \left(\frac{1 -\epsilon}{z}\right)^2}$ & Branch points at $\pm ( 1 - \epsilon )$ \\
        $f_M(z) = c \tan\left(z^{-4}\right)$ & Meromorphic with essential singularity at $0$ \\
    \end{tabular}
    \caption{Test functions used in numerical experiments on the unit circle. The scaling constant $c$ is always chosen so that $\infnorm{f^{(m)}} = 1$ over the domain, where $m$ is the derivative order being tested.}
    \label{tab:testfun-circle}
\end{table}

For each test function $f$, we construct a rational approximant $r$ using the continuum AAA algorithm~\cite{driscollAAARational2024} or continuum greedy TCF~\cite{DriscollGreedyThiele2025pre} to a target accuracy of 1000 times machine epsilon. We then construct a test set $T_f$ consisting of the union of the following sets:
\begin{itemize}
    \item 5000 uniformly spaced points in the domain,
    \item the interpolation nodes of the rational approximant $r$, and
    \item the points $\pm 10^{-16}, \pm 10^{-15.9}, \pm 10^{-15.8}, \ldots, \pm 10^{-3}$ away from each node (in an angular sense for the unit circle), inside the domain.
\end{itemize}
In our results, we report the maximum over $T_f$ of three types of errors:
\begin{itemize}
    \item The absolute error $|f^{(m)}(z) - r^{(m)}(z)|$ using our new algorithms for either representation, called \emph{Formula} error in the tables;
    \item The same quantity computed using extended precision arithmetic (\texttt{Double64} from the DoubleFloats.jl package), called  \emph{Formula-EP} error; and
    \item For the barycentric representation resulting from AAA, we compute the errors for the Schneider--Werner formula~\eqref{eq:sw1} in double and extended precision, called \emph{SW} and \emph{SW-EP} errors, respectively. 
\end{itemize}
We expect to see a difference between double- and extended-precision errors only if an evaluation formula is numerically unstable. 

\begin{table}
    \centering
 \begin{tabular}{l|rrrr|rr}
& \multicolumn{4}{c|}{Barycentric} & \multicolumn{2}{c}{TCF} \\
Function & Formula & Formula-EP & SW & SW-EP & Formula & Formula-EP \\
\hline
\multicolumn{7}{l}{\emph{First derivative, $m=1$}} \\[0.5ex]
$f_E$ & 1.461e-13 & 1.384e-13 & 4.648e-01 & 1.384e-13 & 9.431e-13 & 9.430e-13 \\
$f_C$ & 1.555e-13 & 1.513e-13 & 2.420e-02 & 1.513e-13 & 1.670e-12 & 1.643e-12 \\
\hline
\multicolumn{7}{l}{\emph{Second derivative, $m=2$}} \\[0.5ex]
$f_E$ & 1.853e-11 & 1.486e-11 & 5.289e+01 & 1.486e-11 & 7.863e-11 & 7.863e-11 \\
$f_C$ & 3.884e-12 & 3.858e-12 & 8.620e-02 & 3.858e-12 & 1.083e-11 & 1.081e-11 \\
\end{tabular}
\caption{Results for the nonparametric test functions from \autoref{tab:testfun-interval} on the interval $[-1,1]$.}
\label{tab:errors-interval}
\end{table}

\begin{table}
    \centering
 \begin{tabular}{l|rrrr|rr}
& \multicolumn{4}{c|}{Barycentric} & \multicolumn{2}{c}{TCF} \\
Function & Formula & Formula-EP & SW & SW-EP & Formula & Formula-EP \\
\hline
\multicolumn{7}{l}{\emph{First derivative, $m=1$}} \\[0.5ex]
$f_E$ & 3.047e-12 & 3.046e-12 & 4.409e-01 & 3.046e-12 & 4.541e-12 & 4.541e-12 \\
$f_O$ & 6.381e-13 & 5.975e-13 & 2.162e-01 & 5.975e-13 & 7.331e-12 & 7.322e-12 \\
$f_M$ & 3.966e-13 & 3.935e-13 & 9.674e-02 & 3.935e-13 & 9.439e-13 & 9.387e-13 \\
\hline
\multicolumn{7}{l}{\emph{Second derivative, $m=2$}} \\[0.5ex]
$f_E$ & 4.135e-11 & 4.133e-11 & 3.044e+00 & 4.133e-11 & 5.883e-11 & 5.883e-11 \\
$f_O$ & 3.981e-12 & 3.883e-12 & 7.793e-01 & 3.883e-12 & 6.025e-11 & 6.025e-11 \\
$f_M$ & 1.688e-12 & 1.687e-12 & 3.022e-01 & 1.687e-12 & 3.127e-12 & 3.105e-12 \\
\end{tabular}
\caption{Results for the nonparametric test functions from \autoref{tab:testfun-circle} on the unit circle.}
    \label{tab:errors-circle}
\end{table}

Tables~\ref{tab:errors-interval} and~\ref{tab:errors-circle} show results on the nonparametric test functions for the interval and circle domains, respectively. These results clearly show the difference in stability between the Schneider--Werner formula and our new formula for the barycentric representation: the formulas achieve identical results in extended precision, but only the new formula produces results close to the EP values. The errors in the TCF representation are essentially the same in double and extended precision, and less than an order of magnitude away from the stable errors in the barycentric representation.

\begin{figure}
    \centering
    \includegraphics[width=1.0\textwidth]{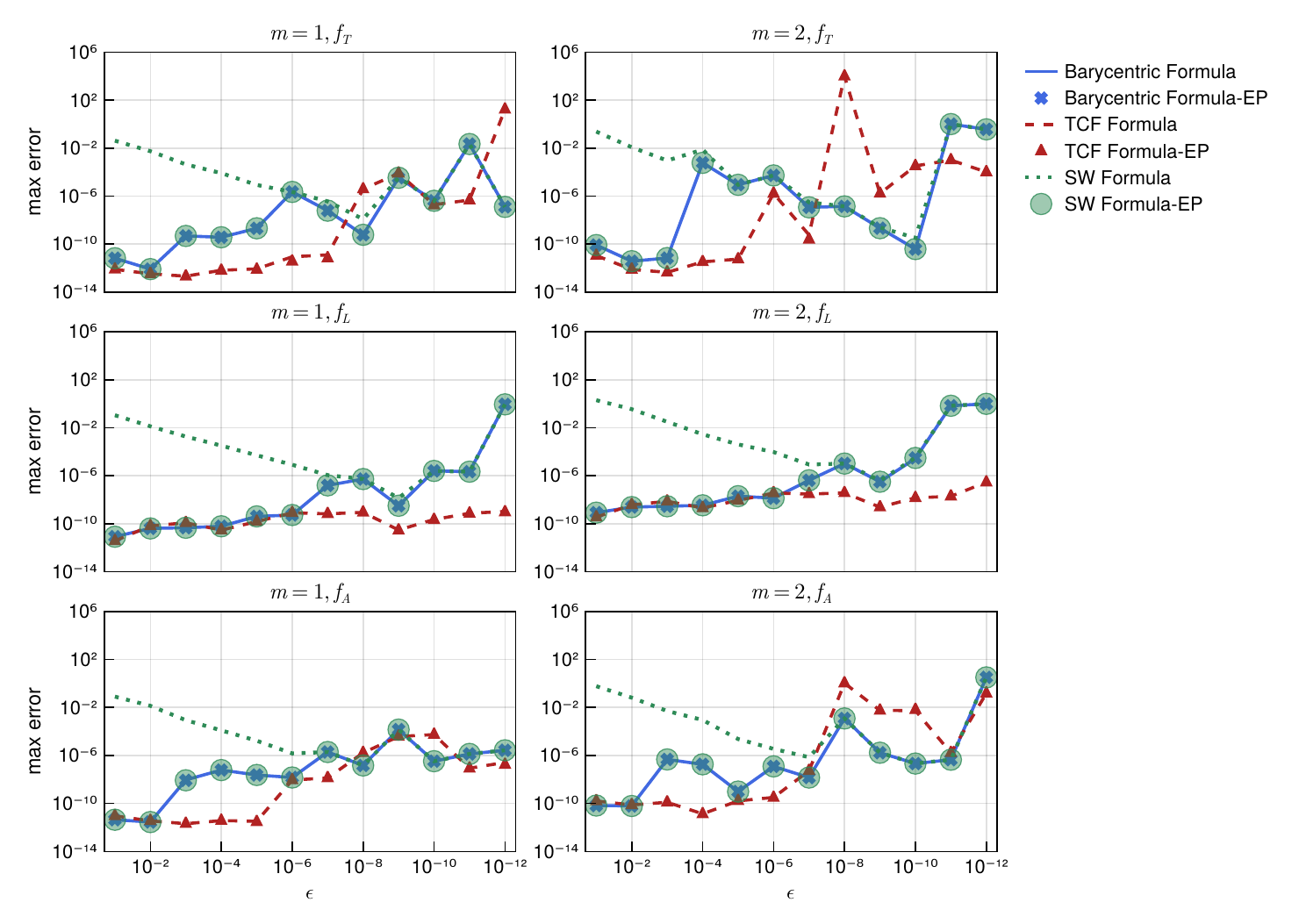}%
    \caption{Errors for evaluating derivatives of parameterized test functions from \autoref{tab:testfun-interval} on the interval $[-1,1]$.}
    \label{fig:errors-interval}
\end{figure}

\begin{figure}
    \centering
    \includegraphics[width=1.0\textwidth]{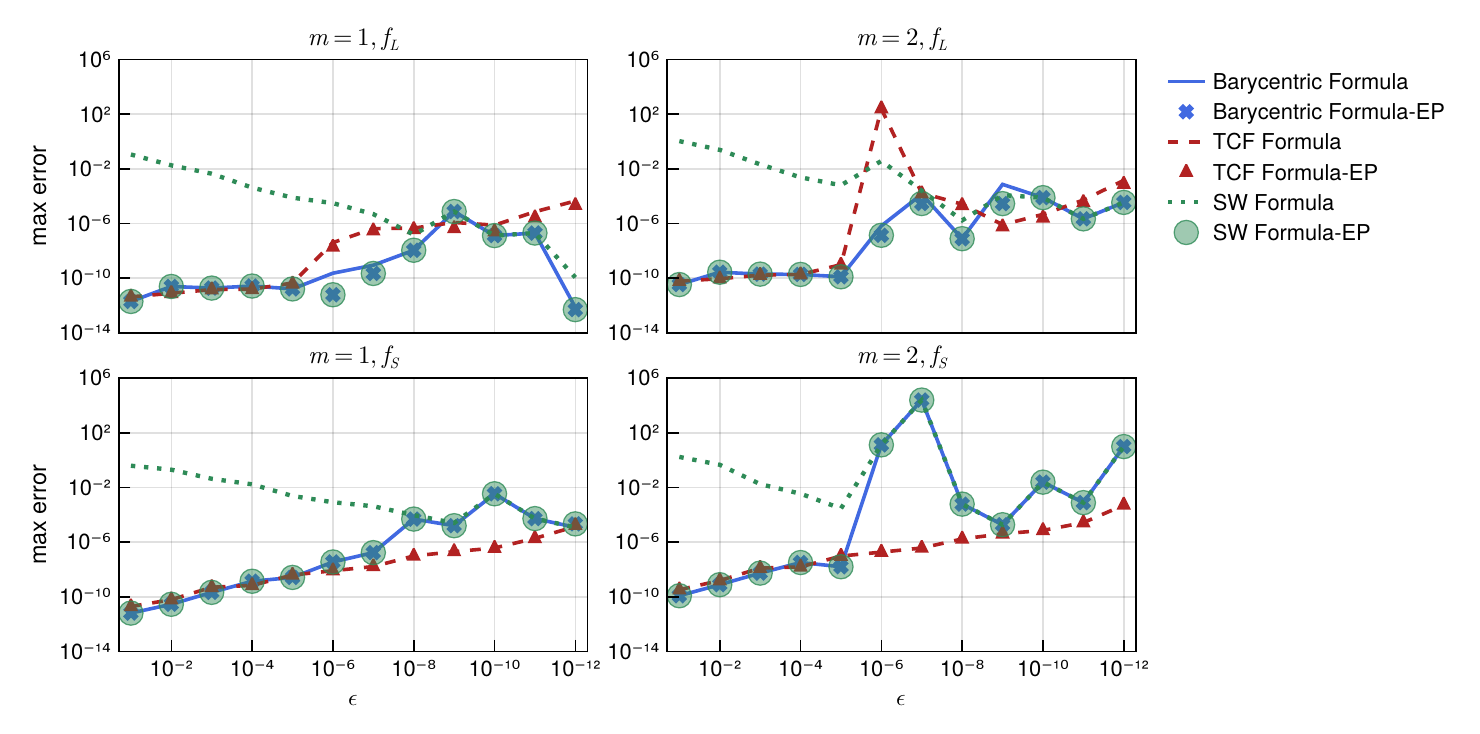}%
    \caption{Errors for evaluating derivatives of parameterized test functions from \autoref{tab:testfun-circle} on the unit circle.}
    \label{fig:errors-circle}
\end{figure}

Figures~\ref{fig:errors-interval} and~\ref{fig:errors-circle} show results for the parameterized test functions on the interval and circle domains, respectively. The results further confirm the numerical stability of the different derivative formula implementations. However, neither the barycentric nor the TCF representation could be considered fully reliable for finding the derivatives globally for $\epsilon \le 10^{-6}$. These are demanding tests: the derivatives of the test functions are a factor $\epsilon^{-m}$ larger than the functions themselves, suggesting large condition numbers to evaluate approximations, and the error is measured in the max norm over a grid much finer than that used for any discrete representation. We conjecture that the occasional failures to get accurate results are due to the imperfections of the continuum algorithms and perhaps to the tendency of both methods to produce rational approximations with poles of small residual extremely close to or on the domain. It may be that a procedure to remove such poles could improve the reliability of both representations.

\section{Conclusions}
\label{sec:conclusions}

We have demonstrated efficient, stable algorithms for evaluating the derivatives of rational functions in both barycentric and Thiele continued fraction representations. The proposed method for the barycentric representation is, to our knowledge, the first globally stable one. While we lack formal proofs of stability, numerical experiments on first and second derivatives confirm it for a variety of test functions and domains. We note that due to the conditioning of the problem, the accuracy obtainable for derivatives of rational approximations can be lower than for the approximations themselves, and for a direct rational approximation of the exact derivative function.

Accounting for the interpolation node nearest the evaluation point $z$ via~\eqref{eq:rderiv} seems to provide sufficient stability for evaluating the first two derivatives of a barycentric representation for all the experiments in \autoref{sec:experiments}. In principle, though, $z$ might be close to two or more nodes when the nodes are tightly clustered. If this proves to be a concern, one might need to resort to the much slower~\eqref{eq:double-sum} for the first derivative and a triple sum for the second derivative, or find more sophisticated approximations, such as multipole expansions. 

\section{Data availability}
No data were used in the research described in the article. Software implementing the (stable) algorithms described is available at \url{https://github.com/complexvariables/RationalFunctionApproximation.jl}.

\section{Disclosure of AI}
Generative AI was used to help write Julia scripts to produce the data for and formatting of error tables and figures.

\backmatter

\bibliography{Approximation}

\end{document}